\documentclass[18pt]{article}
\usepackage[mathscr]{eucal}
\usepackage[english]{babel}
\usepackage{amsmath}
\usepackage{amsfonts}
\usepackage{amsthm}
\usepackage{amssymb}

\usepackage{color}
\usepackage{pstricks, pst-node}

\usepackage{pstricks-add}
\usepackage{graphics,graphicx,graphpap}
\usepackage{xypic}
\usepackage{graphicx}
\usepackage[all,cmtip]{xy}
\usepackage{fancyhdr}
\author{Jaime Castro P\'erez\footnote{jcastrop@itesm.mx,\;corresponding\;author} \\ \textit{Escuela de Ingenier\'ia y Ciencias, Instituto Tecnol\'ologico y de} \\ \textit{Estudios Superiores de Monterrey} \\ \textit{Calle del Puente 222, Tlalpan, 14380, M\'exico D.F., M\'exico.} \\ Mauricio Medina B\'arcenas\footnote{mmedina@matem.unam.mx} \; Jos\'e R\'ios Montes\footnote{jrios@matem.unam.mx} \\  Angel Zald\'ivar\footnote{zaldivar@matem.unam.mx}\\ \textit{Instituto de Matem\'aticas, Universidad Nacional} \\ \textit{Aut\'onoma de M\'exico} \\ \textit{Area de la Investigaci\'on Cient\'ifica, Circuito Exterior, C.U.,} \\ \textit{04510, M\'exico D.F., M\'exico.}}
\title{On Semiprime Goldie Modules}

\begin{document}
\maketitle

\begin{abstract}
For an $R$-module $M$, projective in $\sigma[M]$ and satisfying ascending chain condition (ACC) on left annihilators, we introduce the concept of Goldie module. We also use the concept of semiprime module defined by Raggi et. al. in \cite{S} to give necessary and sufficient conditions for  an $R$-module $M$, to be a semiprime Goldie module. This theorem is a generalization of Goldie's theorem for semiprime left Goldie rings. Moreover, we prove that $M$ is a semiprime (prime) Goldie 
module if and only if the ring $S=End_R(M)$ is a semiprime (prime) right Goldie ring. Also, we study the case when $M$ is a duo module. 
\end{abstract}



\section*{Introduction}

Goldie's Theorem states that a ring $R$ has a semisimple artinian classical left quotient ring if and only if $R$ is a semiprime ring with finite uniform dimension and satisfies ACC on left annihilators. Wisbauer proves in (\cite{W}, Theorem 11.6) a version of Goldie's Theorem in terms of modules. For a retractable $R$-module $M$ with $S=End_R(M)$ the following conditions are equivalent: $\textit{1}$. $M$ is non $M$-singular with finite uniform dimension and $S$ is semiprime, $\textit{2}$. $M$ is non $M$-singular with finite uniform dimension and for every $N\leq_e{M}$ there exists a monomorphism $M\rightarrow{N}$, $\textit{3}$. $End_R(\widehat{M})$ is semisimple left artinian and it is the classical left quotient ring of $S$, here $\widehat{M}$ denotes the $M$-injective hull of $M$. Also, in \cite{H} the authors study when the endomorphism ring of a semiprojective module is a semiprime Goldie ring.

In this paper we give another generalization of Goldie's Theorem. For this, we use the product of submodules of a module $M$ defined in \cite{Bic} to say when a module is a semiprime module. This product extends the product of left ideals of a ring $R$, so $R$ is a semiprime module (over itself) if and only if $R$ is a semiprime ring in the usual sense.

In order to have a definition of Goldie Module such that it extends the classical definition of left Goldie ring, we introduce what ascending chain condition on left annihilators means on a module. A left annihilator in $M$ is a submodule of the form $\mathcal{A}_X=\bigcap_{f\in{X}}{Ker(f)}$ for some $X\subseteq{End_R(M)}$. This definition with $R=M$ is the usual concept of left annihilator. 

The main concept of this work is that an $R$-module $M$ is a Goldie module if $M$ satisfies ACC on left annihilators and has finite uniform dimension. We prove some characterizations of semiprime Goldie modules (Theorem \ref{110}, Theorem \ref{144} and Corollary \ref{145}) which generalize the Goldie's Theorem and extends the Theorem 11.6 of \cite{W} and corollary 2.7 of \cite{H}.

We organize this paper in three sections. Section 1 proves several results for semiprime modules. We also generalize Theorem 10.24 of \cite{LN} to semiprime artinian modules.

    In section 2 we introduce the concept of Goldie modules. We prove the main Theorem of this paper and a characterization of semiprime Goldie modules. We also obtain some examples of Goldie modules. We also prove that if $M$ has finitely many minimal prime submodules $P_1$,...,$P_t$ in $M$ such that $M/P_i$ $(1\leq{i}\leq{t})$ has finite uniform dimension, then $M$ is Goldie module if and only if each $M/P_i$ is Goldie module for $(1\leq{i}\leq{t})$. We also give a description of the submodule $\mathcal{Z}(N)$ with $N\in\sigma[M]$.

In the last section we apply the previous results to duo modules which extend results for commutative rings. In \cite{DM} the authors say that they do not know a duo module with a quotient not duo, in this section we show an example.

Throughout this paper $R$ will be an associative ring with unit and $R$-Mod will denote the category of unitary left $R$-modules. A submodule $N$ of an $R$-module $M$ is denoted by $N\leq{M}$. If $N$ is a proper submodule we write $N<M$. We use $N\leq_e{M}$ for an essential submodule. Let $M$ and $X$ be $R$-modules. $X$ is said to be $M$-generated if there exists an epimorphism from a direct sum of copies of $M$ onto $X$. Every $R$-module $X$ has a largest $M$-generated submodule called the trace of $M$ in $X$, defined by $tr^M(X)=\sum\{f(M)|f:M\to X\}$. The category $\sigma[M]$ is defined as the smallest full subcategory of $R$-Mod containing all $R$-modules $X$ which are isomorphic to a submodule of an $M$-generated module. 

A module $N\in\sigma[M]$ is called singular in $\sigma[M]$ or $M$-singular, if there is an exact sequence in $\sigma[M]$, $0\rightarrow{K}\rightarrow{L}\rightarrow{N}\rightarrow{0}$ with $K\leq_e{L}$. The class 
$\mathcal{S}$ of all $M$-singular modules in $\sigma[M]$ is closed under submodules, quotients and direct sums. Therefore, any $L\in\sigma[M]$ has a largest $M$-singular submodule 
\[\mathcal{Z}(L)=\sum\{f(N)|N\in\mathcal{S}\;\rm{and}\;f\in{Hom_R(N,L)}\}\]
$L$ is called non $M$-singular if $\mathcal{Z}(L)=0$. 

Let $M$ be an $R$-module. In \cite{B} the annihilator in $M$ of a class $\mathcal{C}$ of modules is defined as $Ann_M(\mathcal{C})=\bigcap_{K\in\Omega}{K}$, where
\[\Omega=\{K\leq{M}|\rm{there}\;\rm{exists}\;W\in\mathcal{C}\;and\;f\in{Hom_R(M,W)}\;\rm{with}\;K=Ker(f)\}\]
Also in \cite{B}, the author defines a product in the following way: Let $N\leq{M}$. For each module $X$, $N\cdot{X}=Ann_M(\mathcal{C})$ where $\mathcal{C} $ is the class of modules $W$ such that $f(N)=0$ for all $f\in{Hom_R(M,W)}$. 

For an $R$-module $M$ and $K,L$ submodules of $M$, in \cite{Bic} the product $K_ML$ is defined by $K_ML=\sum\{f(K)|f\in{Hom_R(M,L)}\}$. Moreover, in \cite{B} it is showed that if $M$ is projective in $\sigma[M]$, and $N\leq{M}$, then $N\cdot{X}=N_MX$ for every module $X$. 

A nonzero $R$-module $M$ is called monoform if for each submodule $N$ of $M$ and each morphism $f:N\rightarrow{M}$, $f$ is either zero or a monomorphism. $M$ has enough monoforms if each nonzero submodule of $M$ contains a monoform submodule.

Let $M$-tors be the frame of all hereditary torsion theories on $\sigma[M]$. For a family $\{M_\alpha\}$ of modules in $\sigma[M]$, let $\chi(\{M_\alpha\})$ the greatest element of $M$-tors for which all $M_\alpha$ are torsion free. Let $\xi(\{M_\alpha\})$ be the least element of $M$-tors for which all $M_\alpha$ are torsion. $\xi(\{M_\alpha\})$ and $\chi(\{M_\alpha\})$ are called the hereditary torsion theory generated by the family $\{M_\alpha\}$ and the hereditary torsion theory cogenerated by the same family. In particular, the greatest and least elements in $M$-tors are denoted by $\chi$ and $\xi$ respectively. If $\tau\in{M-tors}$, let $\mathbb{T}_\tau$, $\mathbb{F}_\tau$ and $t_\tau$ denote the torsion class, the torsion free class and the preradical associated to $\tau$, respectively. For details about concepts and terminology concerning torsion theories in $\sigma[M]$, see \cite{Foun} and \cite{W}. 
 
\section{Semiprime Modules}

\newtheorem{11}{Definition}[section]
\begin{11}\label{I11}
\normalfont
Let $M\in{R-Mod}$ and $K$, $L$ submodules of $M$. Put $K_ML=\sum\{f(K)|f\in{Hom_R(M,L)}\}$. For the properties of this product see \cite{P} Proposition  1.3. 
\end{11}

\newtheorem{I2}[11]{Definition}
\begin{I2}
\normalfont
Let $M\in{R-Mod}$. We say a fully invariant submodule $N\leq{M}$ is a \emph{prime submodule} in $M$ if for any fully invariant submodules $K,L\leq{M}$ such that $K_ML\leq{N}$, then $K\leq{N}$ or $L\leq{N}$. We say $M$ is a \emph{prime module}  if $0$ is a prime submodule.
\end{I2}

\newtheorem{31}[11]{Proposition}
\begin{31}\label{121}
Let $M$ be projective in $\sigma[M]$ and $P$ a fully invariant submodule of $M$. The following conditions are equivalent:
\begin{enumerate}		
	\item $P$ is prime in $M$.
	\item For any submodules $K$, $L$ of $M$ containing $P$ and such that $K_ML\leq{P}$, then $K=P$ or $L=P$.
\end{enumerate}
\end{31}

\begin{proof}
$\textit{1}\Rightarrow\textit{2}:$ By Proposition 1.11 of \cite{P}. 

$\textit{2}\Rightarrow\textit{1}:$ Suppose that $K$, $L$ are submodules of $M$ such that $K_ML\leq{P}$. 

We claim that $K_M(L+P)\leq{P}$. Since $K_ML\leq{L\cap{P}}$, by Proposition 5.5 of \cite{B} $K_M({L}/{L\cap{P}})=0$ so $K_M({L+P}/{P})=0$. Thus $K_M(L+P)\leq{P}$.

On the other hand, 
\[(K+P)_M(L+P)=K_M(L+P)+P_M(L+P)\leq{P}\]
because $P$ is fully invariant in $M$.

Then, by hypothesis $K+P={P}$ or $L+P={P}$, hence $K\leq{P}$ or $L\leq{P}$.
\end{proof}

\newtheorem{I3}[11]{Definition}
\begin{I3}
\normalfont
We say a fully invariant submodule $N\leq{M}$ is a \emph{semiprime submodule}  in $M$ if for any fully invariant submodule $K\leq{M}$ such that $K_MK\leq{N}$, then $K\leq{N}$. We said $M$ is a \emph{semiprime module}  if $0$ is a semiprime submodule.
\end{I3}    

\newtheorem{53}[11]{Lemma}
\begin{53}\label{143}
Let $M$ be projective in $\sigma[M]$ and $N$ a fully invariant submodule of $M$. The following conditions are equivalent:
\begin{enumerate}
	\item $N$ is semiprime in $M$.
	\item For any submodule $K$ of $M$, $K_MK\leq{N}$ implies $K\leq{M}$.
	\item For any submodule $K\leq{M}$ containing $N$ such that $K_MK\leq{N}$, then $K=N$.
\end{enumerate}
\end{53}
\begin{proof}
$\textit{1}\Rightarrow\textit{2}:$ Let $K\leq{M}$ such that $K_MK\leq{N}$. Consider the submodule $K_MM$ of $M$. This is the minimal fully invariant submodule of $M$ which contains $K$ and $K_MX=(K_MM)_MX$ for every module $X$. Hence by Proposition 1.3 of \cite{P}  we have that
\[K_MK=(K_MM)_MK\leq{((K_MM)_MK)_MM)}\leq{N_MM}\]
Since $N$ is fully invariant submodule of $M$ then $N_MM=N$ and by Proposition  5.5 of \cite{B} $(K_MM)_M(K_MM)=((K_MM)_MK)_MM)\leq{N}$. Since $N$ is semiprime in $M$, $K_MM\leq{N}$. Hence $K\leq{N}$.

$\textit{2}\Rightarrow\textit{1}:$ By definition.

$\textit{1}\Leftrightarrow\textit{3}:$ Similar to the proof of Proposition \ref{121}.
\end{proof}

In Remark \ref{aa} below, we give an example where the associativity of the product $(\cdot)_M(\cdot)$ is not true in general.

\newtheorem{I4}[11]{Definition}
\begin{I4}
\normalfont
Let $M\in{R-Mod}$ and $N$ a fully invariant submodule of $M$. We define the \emph{powers}  of $N$ as:
\begin{enumerate}
	\item $N^0=0$
	\item $N^1=N$
	\item $N^m=N_MN^{m-1}$
\end{enumerate}
\end{I4}

\newtheorem{34}[11]{Lemma}
\begin{34}\label{124}
Let $M$ be projective in $\sigma[M]$ and $N$ semiprime in $M$. Let $J$ be a fully invariant submodule of $M$ such that $J^n\leq{N}$ then $J\leq{N}$.
\end{34}

\begin{proof}
By induction on $n$. If $n=1$ the result is clear.

Suppose $n>1$ and the Proposition is valid for $n-1$. We have that $2n-2\geq{n}$ then
\[J^{2n-2}\leq{N}\]
so
\[(J^{n-1})^2={J^{n-1}}_M{J^{n-1}}\leq{N}\]
since $N$ is semiprime $J^{n-1}\leq{N}$ then $J\leq{N}$.
\end{proof}

\newtheorem{18}[11]{Proposition}
\begin{18}\label{108a}
Let $S:=End_R(M)$ and assume $M$ generates all its submodules. If $N$ is a fully invariant submodule of $M$ such that $Hom_R(M,N)$ is a prime (semiprime) ideal of $S$, then $N$ is prime (semiprime) in $M$.
\end{18}
\begin{proof}
Let $K$ and $L$ be fully invariant submodules of $M$ such that $K_ML\leq{N}$. Put $I=Hom_R(M,L)$ and $J=Hom_R(M,K)$. Let $m\in{M}$ and $\sum{f_ig_i}\in{IJ}$. Since $g_i\in{J}$ and $g_i(m)\in{K}$ then $\sum{f_i(g_i(m))}\in{K_ML}\leq{N}$. Hence $IJ\leq{Hom_R(M,N)}$. Since $Hom_R(M,N)$ is prime (semiprime) in $S$, then $I\leq{Hom_R(M,N)}$ or $J\leq{Hom_R(M,N)}$. Hence $tr^M(L):=Hom(M,L)M\leq{N}$ or $tr^M(K)\leq{N}$ and since $M$ generates all its submodules then $L\leq{N}$ or $K\leq{N}$. Thus $N$ is a prime (semiprime) submodule. 
\end{proof}

Next definition aper in \cite{Kh}

\newtheorem{18a}[11]{Definition}
\begin{18a}
A module $M$ is \emph{retractable} if $Hom_R(M,N)\neq{0}$ for all $0\neq{N}\leq{M}$
\end{18a}

\newtheorem{I5}[11]{Corollary}
\begin{I5}\label{108}
Let $S:=End_R(M)$ with $M$ retractable. If $S$ is a prime (semiprime) ring then $M$ is prime (semiprime).
\end{I5}
\begin{proof}
Let $K$ and $L$ be fully invariant submodules of $M$ such that $K_ML=0$. Since $Hom_R(M,0)$ is a prime (semiprime) ideal of $S$ then  by the proof of \ref{108a}, $tr^M(K)=0$ o $tr^M(L)=0$. Since $M$ is retractable, $K=0$ or $L=0$. Hence  $0$ is prime (semiprime) in $M$. Thus $M$ is prime (semiprime).
\end{proof}

\newtheorem{111}[11]{Proposition}
\begin{111}\label{123}
Let $M$ be projective in $\sigma[M]$ and $N$ a proper fully invariant submodule of $M$. The following conditions are equivalent:
\begin{enumerate}
	\item $N$ is semiprime in $M$.
	\item If $m\in{M}$ is such that ${Rm}_M{Rm}\leq{N}$, then $m\in{N}$.
	\item $N$ is an intersection of prime submodules.
\end{enumerate}
\end{111}

\begin{proof}
$\textit{1}\Rightarrow\textit{2}:$ By Lemma \ref{143}. 

$\textit{2}\Rightarrow\textit{3}:$ Since $N$ is proper in $M$, let $0\neq{m_0}\in{M\setminus{N}}$. Then ${Rm_0}_M{Rm_0}\nleq{N}$. Now, let $0\neq{m_1}\in{{Rm_0}_M{Rm_0}}$ but $m_1\notin{N}$ Then ${Rm_1}_M{Rm_1}\nleq{N}$ and ${Rm_1}_M{Rm_1}\leq{{Rm_0}_M{Rm_0}}$. We obtain a sequence of non-zero elements of $M$, $\{m_0,m_1,...\}$ such that $m_i\notin{N}$ for all $i$ and ${Rm_{i+1}}_M{Rm_{i+1}}\leq{{Rm_i}_M{Rm_i}}$.

By Zorn's Lemma there exists a fully invariant submodule $P$ of $M$ with $N\leq{P}$, maximal with the property that $m_i\notin{P}$ for all $i$ . 

We claim $P$ is a prime submodule. Let $K$ and $L$ submodules of $M$ containing $P$. Since $P\leq{K}$ and $P\leq{L}$, then there exists $m_i$ and $m_j$ such that $m_i\in{K}$ and $m_j\in{L}$. Suppose $i\leq{j}$, then ${Rm_i}_M{Rm_i}\leq{K}$ and by construction $m_j\in{Rm_i}_M{Rm_i}$ and thus $m_j\in{K}$. If we put $k=max\{i,j\}$, then $m_k\in{K}$ and $m_k\in{L}$. Hence, ${Rm_k}_M{Rm_k}\leq{K_ML}$, and so $K_ML\nleq{P}$. By Proposition  \ref{121}, $P$ is prime in $M$. 

$\textit{3}\Rightarrow\textit{1}:$ It is clear. 
\end{proof}

\newtheorem{26c}[11]{Proposition}
\begin{26c}\label{116c}
Let $0\neq{M}$ be a semiprime module and projective in $\sigma[M]$. Then $M$ has minimal prime submodules in $M$.
\end{26c}

\begin{proof}
By the proof Proposition of \ref{123}, $M$ has prime submodules. Let $P\leq{M}$ be a prime submodule. Consider $\Gamma=\{Q\leq{P}|Q\;is\;prime\}$. This family is not empty because $P\in\Gamma$. Let $\mathcal{C}=\{Q_i\}$ be a descending chain in $\Gamma$. Let $N,K\leq{M}$ be fully invariant submodules of $M$ such that $N_MK\leq\bigcap\mathcal{C}$. Suppose that $N\nleq\bigcap\mathcal{C}$. Then there exists $Q_j$ such that $N\nleq{Q_j}$ and $N\nleq{Q_l}$ for all $Q_l\leq{Q_j}$. Therefore $K\leq{Q_l}$ for all $Q_l\leq{Q_j}$, and since $\mathcal{C}$ is a chain then $K\leq\bigcap\mathcal{C}$. Therefore $\bigcap\mathcal{C}\in\Gamma$. By Zorn's Lemma $\Gamma$ has minimal elements.
\end{proof}
 
\newtheorem{26d}[11]{Remark}
\begin{26d}
\normalfont
Notice that if $M$ is projective in $\sigma[M]$ and $M$ has prime submodules in $M$, then $M$ has minimal prime submodules.
\end{26d}

\newtheorem{26a}[11]{Corollary}
\begin{26a}\label{116a}
Let $0\neq{M}$ be a semiprime module and projective in $\sigma[M]$. Then 
\[0=\bigcap\{P\leq{M}|P\;is\;a\;minimal\;prime\;in\;M\}.\]
\end{26a}

\begin{proof}
Let $x\in\bigcap\{P\leq{M}|P\;is\;a\;minimal\;prime\;in\;M\}$ and $Q\leq{M}$ be a prime submodule in $M$. By Proposition  \ref{116c} there exists a minimal prime submodule $P$ such that $P\leq{Q}$ then $x\in{Q}$ and $x$ is in the intersection of all primes in $M$. By Proposition  \ref{123}, $x=0$.
\end{proof}

\newtheorem{36}[11]{Lemma}
\begin{36}\label{126}
Let $M\in{R-Mod}$ and $N$ a minimal submodule of $M$. Then $N^2=0$ or $N$ is a direct summand of $M$.
\end{36}

\begin{proof}
Suppose that $N_MN\neq{0}$. Then there exists $f:M\rightarrow{N}$ such that $f(N)\neq{0}$. Since $0\neq{f(M)}\leq{N}$ and $N$ is a minimal submodule, $f(M)=N$. On the other hand, $Ker(f)\cap{N}\leq{N}$, since $f(N)\neq{0}$ then $Ker(f)\cap{N}=0$. We have that $({M}/{Ker(f)})\cong{N}$ and since $N$ is a minimal submodule $Ker(f)$, then is a maximal submodule of $M$. Thus $Ker(f)\oplus{N}=M$.
\end{proof}  

\newtheorem{37}[11]{Corollary}
\begin{37}\label{127}
Let $M$ be a retractable module. If $N$ is a minimal submodule in a semiprime module $M$, then $N$ is a direct summand.
\end{37}

\begin{proof}
Since $M$ is semiprime, $N_MN\neq{0}$. 
\end{proof}

\newtheorem{38}[11]{Theorem}
\begin{38}\label{128}
The following conditions are equivalent for a retractable $R$-module $M$:
\begin{enumerate}
	\item $M$ is semisimple and left artinian.
	\item $M$ is semiprime and left artinian.
	\item $M$ is semiprime and satisfies DCC on cyclic submodules and direct summands.
\end{enumerate}
\end{38}

\begin{proof}
$\textit{1}\Rightarrow\textit{2}:$ If $M$ is semisimple then it is semiprime. 

$\textit{2}\Rightarrow\textit{3}:$ Since $M$ is left artinian, then it satisfies DCC on cyclic submodules and direct summands. 

$\textit{3}\Rightarrow\textit{1}:$ Since $M$ satisfies DCC on cyclic submodules, there exists $K_1$ a minimal submodule of $M$. By Corollary \ref{127}, $M=K_1\oplus{L_1}$. Now there exists $K_2$ a minimal submodule of $L_1$ and $L_1=K_2\oplus{L_2}$. With this process we obtain a descending chain of direct summands, which by hypothesis it is finite $L_1\supseteq{L_2}\supseteq{L_3}\supseteq...\supseteq{L_m}$. Since $L_m$ is simple and $M=K_1\oplus{K_2}\oplus...\oplus{K_m}\oplus{L_m}$, then $M$ is semisimple.  

Now, if $M$ is semisimple and satisfies DCC on direct summands then $M$ is artinian.
\end{proof}

\newtheorem{21}[11]{Definition}
\begin{21}\label{11}
\normalfont
Let $M\in{R-Mod}$ and $N\leq{M}$. We say $N$ is an \emph{annihilator submodule}  if $N=Ann_M(K)$ for some $0\neq{K}\leq{M}$.
\end{21}

\newtheorem{22a}[11]{Lemma}
\begin{22a}\label{112a}
Let $M$ be semiprime and projective in $\sigma[M]$. Let $N,L\leq{M}$. If $L_MN=0$, then $N_ML=0$ and $L\cap{N}=0$.
\end{22a}
\begin{proof}
Since $L_MN=0$, then 
\[0=N_M(L_MN)_ML=(N_ML)_M(N_ML).\]
Hence $N_ML=0$ . 

Now, since $L\cap{N}\leq{L}$ and $L\cap{N}\leq{N}$, then
\[(L\cap{N})_M(L\cap{N})\leq{L_MN}=0.\]
Thus $L\cap{N}=0$
\end{proof}

\newtheorem{22}[11]{Corollary}
\begin{22}\label{112}
Let $M$ be semiprime and projective in $\sigma[M]$. If $N\leq{M}$, then $N_MAnn_M(N)=0$.
\end{22}

\newtheorem{23}[11]{Proposition}
\begin{23}\label{113}
Let $M$ be semiprime and projective in $\sigma[M]$ and $N\leq{M}$. Then $N$ is an annihilator submodule if and only if $N=Ann_M(Ann_M(N))$
\end{23}
\begin{proof}
$\Rightarrow:$ By Lemma \ref{112} $N\leq{Ann_M(Ann_M(N))}$. There is $K\leq{M}$ such that $N=Ann_M(K)$, hence 
\[K_MN=K_MAnn_M(K)=0\]
and thus $K\leq{Ann_M(N)}$. Therefore,
\[Ann_M(Ann_M(N))\leq{Ann_M(K)}=N\]
It follows that $N=Ann_M(Ann_M(N))$. 

$\Leftarrow:$ By definition of annihilator submodule.
\end{proof}

\newtheorem{24}[11]{Proposition}
\begin{24}\label{114}
Let $M$ be semiprime and $N\leq{M}$. Then, $Ann_M(N)$ is the unique pseudocomplement fully invariant of $N$. Moreover, $N\bigoplus{Ann_M(N)}$ intersects all fully invariant submodules of $M$.
\end{24}
\begin{proof}
Let $L\leq{M}$ be a fully invariant pseudocomplement of $N$ in $M$. Then 
\[L_MN\leq{L\cap{N}}=0\]
Thus $L\leq{Ann_M(N)}$. Observe that
\[(Ann_M(N)\cap{N})_M(Ann_M(N)\cap{N})\leq(Ann_M(N)\cap{N})_MN=0\]
Since $M$ is semiprime, $Ann_M(N)\cap{N}=0$.
Thus $L=Ann_M(N)$.
\end{proof}

\newtheorem{26b}[11]{Lemma}
\begin{26b}\label{116b}
Let $M$ be a semiprime module and $N\leq{M}$. Let $S$ be the set of all minimal prime submodules of $M$ which do not contain $N$. Then $Ann_M(N)=\bigcap\{P|P\in{S}\}$.
\end{26b}

\begin{proof}
Put $K=\bigcap\{P|P\in{S}\}$. Any element in $K\cap{N}$ is in the intersection of all minimal prime submodules of $M$ which is zero. Then $K\cap{N}=0$. Since $K$ is fully invariant in $M$, $K_MN\leq{K\cap{N}}=0$. Thus, $K\leq{Ann_M(N)}$. Now, let $P\in{S}$. Since $Ann_M(N)_MN=0\leq{P}$ and $N\nleq{P}$, then $Ann_M(N)\leq{K}$.
\end{proof}

\newtheorem{59}[11]{Lemma}
\begin{59}\label{148}
Let $M$ be projective in $\sigma[M]$. If $M$ is semiprime then $M$ is retractable.
\end{59}

\begin{proof}
Let $N\leq{M}$ and suppose $Hom_R(M,N)=0$. Then $Ann_M(N)=M$. So $M_MN=0$ but $N_MN\subseteq{M_MN}=0$. Since $M$ is semprime then $N=0$ by Lemma \ref{143}.
\end{proof}

%In the rest of this section assume that $M\in{R-Mod}$ is projective in $\sigma[M]$.

\newtheorem{25}[11]{Proposition}
\begin{25}\label{115}
Let $M$ be projective in $\sigma[M]$ and semiprime. The following conditions are equivalent for $N\leq M$:
\begin{enumerate}
	\item $N$ is a maximal annihilator submodule.
	\item $N$ is an annihilator submodule and is a minimal prime submodule.
	\item $N$ is prime in $M$ and $N$ is an annihilator submodule.
\end{enumerate} 
\end{25}
\begin{proof}
$\textit{1}\Rightarrow\textit{2}:$ Let $K\leq{M}$ such that $N=Ann_M(K)$. Let $L,H\leq{M}$ be fully invariant submodules of $M$ such that $L_MH\leq{N}$. Assume $H\nleq{N}$. Then $0\neq{H_MK}$. Hence $Ann_M(K)\leq{Ann_M(H_MK)}$, but since $Ann_M(K)$ is a maximal annihilator submodule, then $Ann_M(K)=Ann_M(H_MK)$. 

As $M$ is projective in $\sigma[M]$, by Proposition 5.5 of \cite{B}, we have that 
\[L_M(H_M(H_MK))=(L_MH)_M(H_MK)\leq{N_M(H_MK)}=0\]
Now, since $H_M(H_MK)\leq{H_MK}$, then
\[Ann_M(K)=Ann_M(H_MK)\leq{Ann_M(H_M(H_MK))}\]
Therefore $Ann_M(H_MK)=Ann_M(H_M(H_MK))$. Thus $L\leq{Ann_M(K)}=N$. 

Now, let $P\leq{M}$ be a prime submodule of $M$ such that $P<{N}$. We have that $N_MK=0\leq{P}$. So $K\leq{P}<{N}$. Hence $K_MK=0$. Thus $K=0$, but $M$ is semiprime, a contradiction. It follows that $N$ is a minimal prime submodule of $M$.

$\textit{2}\Rightarrow\textit{3}:$ By hypothesis. 

$\textit{3}\Rightarrow\textit{1}:$ Suppose $N<{K}$ with $K$ an annihilator submodule. Then
\[Ann_M(K)_MK=0\leq{N}\]
Since $N$ is prime in $M$, then $Ann_M(K)\leq{N}<{K}$. By Proposition  \ref{114} $Ann_M(K)\cap{K}=0$, hence $Ann_M(K)=0$. Since $K$ is an annihilator submodule, by Proposition  \ref{113}, $K=Ann_M(Ann_M(K))=Ann_M(0)=M$.
\end{proof}

\newtheorem{1a}[11]{Remark}  
\begin{1a}\label{aa}
\normalfont
Following the notation of Example 1.12 of \cite{P}   
Let $R=\mathbb{Z}_2\rtimes(\mathbb{Z}_2\oplus\mathbb{Z}_2)$. This ring has only one maximal ideal $I$ and it has three simple ideals: $J_1$, $J_2$, $J_3$, which are isomorphic. Then, the lattice of ideals of $R$ has the form 
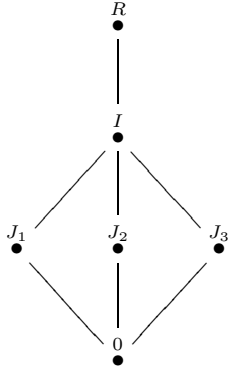
\[\xymatrix{ & \stackrel{R}{\bullet}\ar@{-}[d] & \\ & \stackrel{I}{\bullet}\ar@{-}[d]\ar@{-}[dl]\ar@{-}[dr] & \\ \stackrel{J_1}{\bullet}\ar@{-}[dr] & \stackrel{J_2}{\bullet}\ar@{-}[d] & \stackrel{J_3}{\bullet}\ar@{-}[dl] \\ & \stackrel{0}{\bullet} &}\] 

Moreover, $R$ is artinian and $R$-Mod has only one simple module up to isomorphism. Let $S$ be a simple module. By Theorem 2.13 of \cite{BP}, the lattice of fully invariant submodules of $E(S)$ has tree maximal submodules $N$, $L$ and $K$, and it has the form

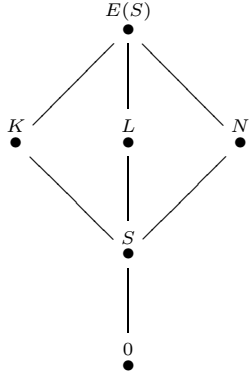
\[\xymatrix{ & \stackrel{E(S)}{\bullet}\ar@{-}[dl]\ar@{-}[d]\ar@{-}[dr] & \\ \stackrel{K}{\bullet}\ar@{-}[dr] & \stackrel{L}{\bullet}\ar@{-}[d] & \stackrel{N}{\bullet}\ar@{-}[dl] \\ & \stackrel{S}{\bullet}\ar@{-}[d] & \\ & \stackrel{0}{\bullet} & }\]

Put $M=E(S)$. Since $K\cap{L}=S$ and $K_ML\leq{K\cap{L}}$, then $K_ML\leq{S}$. On the other hand consider the composition

\[\xymatrix{M\ar[r]^-{\pi} & M/N\ar[r]^-{\cong} & S\ar[r]^-{i} & L}\]

$f=i\circ\pi$ where $\pi$ is the natural projection and $i$ is the inclusion. Then, $f(K)=S$ and $S\leq{K_ML}$. Thus, $K_ML=S$. Notice that $K_ML\leq{N}$ but $K\nleq{N}$ and $L\nleq{N}$. Hence $N$ is not prime in $M$. Analogously, we prove that neither $K$ nor $L$ are prime in $M$. We also note that $K_MK=S$. Moreover, $\pi(K)=S$, so $K_MS=S$. In the same way $L_MS=S$ and $N_MS=S$ 

Let $g:M\rightarrow{K}$ be a non zero morphism. If $Ker(g)\cap{S}=0$ then $g$ is a monomorphism, a contradiction. So $Ker(g)\cap{S}=S$. Thus $S_MK=0$ and $Ann_M(K)=S$. Analogously $Ann_M(L)=S=Ann_M(N)=Ann_M(S)$. Since $S_MS\leq{S_MK}$, $S_MS=0$. Thus $M$ is not semiprime. Hence, $S$ is a maximal annihilator submodule of $M$ which is not prime because $K_MK=S$.
With this we can see that associativity is  not true in general, because $L_M(K_MS)=L_MS=S$ and $(L_MK)_MS=S_MS=0$.
Notice that, in this example $Hom_R(M,H)\neq{0}$ for all $H\in\sigma[M]$ in particular $M$ is retractable, but $M$ is not projective in $\sigma[M]$.
\end{1a}

\newtheorem{25b}[11]{Proposition}
\begin{25b}\label{115b}
Let $M$ be projective in $\sigma[M]$ and semiprime. For $N\leq M$, if $N=Ann_M(U)$ with $U\leq{M}$ a uniform submodule, then $N$ is a maximal annihilator in $M$.
\end{25b}

\begin{proof}
Suppose that $N<{K}$ with $K$ an annihilator submodule in $M$. Since $N=Ann_M(U)$ by Proposition  \ref{114}, $K\cap{U}\neq{0}$. By hypothesis $U$ is uniform and thus $K\cap{U}\leq_{e}U$. Then
\[(K\cap{U})\oplus{Ann_M(U)}\leq_{e}U\oplus{Ann_M(U)}\]
Now, notice that if $L\leq_{FI}M$, by Proposition  \ref{114} $(U\oplus{Ann_M(U)})\cap{L}\neq{0}$. So $((K\cap{U})\oplus{Ann_M(U)})\cap{L}\neq{0}$. Therefore, $K\cap{L}\neq{0}$ and $K$ intersects all fully invariant submodules of $M$ . Since $K\cap{Ann_M(K)}=0$ and $Ann_M(K)\leq_{FI}M$, then $Ann_M(K)=0$. Thus, $K=Ann_M(Ann_M(K))=Ann_M(0)=M$.
\end{proof}

\newtheorem{27}[11]{Proposition}
\begin{27}\label{117}
Let $M$ be projective in $\sigma[M]$ and semiprime with finite uniform dimension. Then:
\begin{enumerate}
	\item $M$ has finitely many minimal prime submodules.
	\item The number of annihilators submodules is finite.
	\item $M$ satisfies ACC on annihilators submodules. 
\end{enumerate}
\end{27}

\begin{proof}
$\textit{1}:$ Let $U_1,..,U_n$ be uniform submodules of $M$ such that $U_1\oplus...\oplus{U_n}\leq_eM$. By Proposition s \ref{115} and \ref{115b}, $P_i:=Ann_M(U_i)$ is a minimal prime submodule of $M$ for each $i$. By Proposition  \ref{114}, $(U_1\oplus...\oplus{U_n})\cap{Ann_M(U_1\oplus...\oplus{U_n})}=0$ and $P_1\cap...\cap{P_n}\leq{Ann_M(U_1\oplus...\oplus{U_n})}=0$. 

Now, if $P$ is a minimal prime submodule of $M$, then 
\[{P_1}_M{P_2}_M...{_M}P_n\leq{P_1\cap...\cap{P_n}}=0\leq{P}\]
Hence, there exists $j$ such that $P_j\leq{P}$, a contradiction. 

$\textit{2}:$ By Lemma \ref{116b}. 

$\textit{3}:$ It is clear by $\textit{2}$.
\end{proof}

\section{Goldie Modules}

The following definition was taken from \cite{E}
\newtheorem{112}{Definition}[section]
\begin{112}
\normalfont
Let $M\in{R-Mod}$. $M$ is \emph{essentially compressible}  if for every essential submodule $N\leq_e{M}$ there exists a monomorphism $M\rightarrow{N}$.
\end{112}

\newtheorem{I6}[112]{Definition}
\begin{I6}
\normalfont
Let $M\in{R-Mod}$. We call a left annihilator in $M$ a submodule 
\[\mathcal{A}_X=\bigcap\{Ker(f)|f\in{X}\}\]
for some $X\subseteq{End_R(M)}$.
\end{I6}

\newtheorem{113}[112]{Definition}
\begin{113}
\normalfont
We say $M$ is a \emph{Goldie module}  if it satisfies ACC on left annihilators and has finite uniform dimension.
\end{113}

\newtheorem{17}[112]{Lemma}
\begin{17}\label{107}
Suppose $M$ is projective in $\sigma[M]$. If $N\in\sigma[M]$ is essentially compressible, then $Ann_M(N)$ is a semiprime submodule of $M$.
\end{17}
\begin{proof}
Let $L\leq{M}$ be a fully invariant submodule of $M$ such that $L_ML\leq{Ann_M(N)}$. Put 
\[\Gamma=\{K\leq{N}|L_MK=0\}\]
Then $\Gamma\neq\emptyset$ and by Zorn's Lemma there exists a maximal independent family $\{K_i\}_I$ in $\Gamma$. Notice that $\bigoplus_I{K_i}\in\Gamma$ because
\[L_M\bigoplus_I{K_i}=\bigoplus_I{L_MK_i}=0\]
Let $0\neq{A}\leq{N}$ be a submodule. Since $(L_ML)_MA=0$ then $L_MA\in\Gamma$. 

If $L_MA=0$ then $A\in\Gamma$ and $A\cap\bigoplus_I{K_i}\neq{0}$ because $\{K_i\}$ is a maximal independent family in $\Gamma$. 

Now,  if $L_MA\neq{0}$ we also have $(L_MA)\cap\bigoplus_I{K_i}\neq{0}$ and $(L_MA)\cap\bigoplus_I{K_i}\leq{A\cap\bigoplus_I{K_i}}$. Thus $\bigoplus_I{K_i}\leq_eN$. 

By hypothesis there exists a monomorphism $\theta:N\rightarrow\bigoplus_I{K_i}$. Then
\[\theta(L_MN)\leq{L_M\bigoplus_I{K_i}}=0\]
and hence $L_MN=0$. Thus $L\leq{Ann_M(N)}$.
\end{proof}

\newtheorem{52}[112]{Proposition}
\begin{52}\label{142}
Let $M$ be projective in $\sigma[M]$. If $N\in\sigma[M]$ is an $M$-singular module, then $Ker(f)\leq_e{M}$ for all $f\in{Hom_R(M,N)}$.
\end{52}
\begin{proof}
Let $f\in{Hom_R(M,N)}$. Since $N$ is $M$-singular, there exists an exact sequence 
\[\xymatrix{0\ar[r] & K\ar[r]^i & L\ar[r]^\pi & N\ar[r] & 0}\]
in $\sigma[M]$ with $K\leq_e{L}$. Since $M$ is projective in $\sigma[M]$, there exists $\hat{f}:M\rightarrow{L}$ such that $\pi\hat{f}=f$: 
\[\xymatrix{ & M\ar[d]^f\ar[dl]_{\hat{f}} & \\ L\ar[r]^\pi & N\ar[r] & 0}\]
As $K\leq_e{L}$, then $\hat{f}^{-1}(K)\leq_e{M}$. Then
\[f(\hat{f}^{-1}(K))=\pi(\hat{f}(\hat{f}^{-1}(K)))\leq\pi(K)=0.\]
Therefore, $\hat{f}^{-1}(K)\leq{Ker(f)}$ and hence $Ker(f)\leq_e{M}$. 
\end{proof}

\newtheorem{35}[112]{Proposition}
\begin{35}\label{125}
Let $M$ be projective in $\sigma[M]$. If $M$ is essentially compressible then $M$ is non $M$-singular.
\end{35}

\begin{proof}
Suppose $\mathcal{Z}(M)\neq{0}$. If $\mathcal{Z}(M)\leq_eM$, then there exists a monomorphism $\theta:M\rightarrow{\mathcal{Z}(M)}$, by Proposition  \ref{142} $Ker\theta\leq_eM$, a contradiction. Therefore $\mathcal{Z}(M)$ has a pseudocomplement $K$ in $M$ and thus $\mathcal{Z}(M)\oplus{K}\leq_eM$. Hence, there exists a monomorphism $\theta:M\rightarrow{\mathcal{Z}(M)\oplus{K}}$. Let $\pi:\mathcal{Z}(M)\oplus{K}\rightarrow{\mathcal{Z}(M)}$  be the canonical projection, then $Ker(\pi\theta)\leq_eM$ and so $Ker(\pi\theta)=\theta^{-1}(Ker\pi)=\theta^{-1}(K)\leq_eM$. But $\mathcal{Z}(M)\cap\theta^{-1}(K)=0$,contradiction. Thus $\mathcal{Z}(M)=0$.
\end{proof} 

\newtheorem{42}[112]{Lemma}
\begin{42}\label{132}
Let $M\in{R-Mod}$ with finite uniform dimension. Then, for every monomorphism $f:M\rightarrow{M}$, $Im(f)\leq_e{M}$.
\end{42}

\begin{proof}
Let $f:M\rightarrow{M}$ be a monomrfism. If the uniform dimension of $M$ is $n$, $(Udim(M)=n)$ and there exists $K\leq{M}$ such that $f(M)\cap{K}=0$, then $Udim(f(M)\oplus{K})=n+1$, a contradiction. 
\end{proof}

%In the rest of this section assume $M\in{R-Mod}$ is projective in $\sigma[M]$. 
 
\newtheorem{20}[112]{Theorem}
\begin{20}\label{110}
Let $M$ be projective in $\sigma[M]$ with finite uniform dimension. The following conditions are equivalent:
\begin{enumerate}
	\item $M$ is semiprime and non $M$-singular
	\item $M$ is semiprime and satisfies ACC on annihilators
	\item Let $N\leq{M}$, then $N\leq_e{M}$ if and only if there exists a monomorphism $f:M\rightarrow{N}$.
\end{enumerate}
\end{20}
\begin{proof}
$\textit{1}\Rightarrow\textit{2}:$ Since $M$ is non $M$-singular and has finite uniform dimension then, by Proposition  3.6 of \cite{K} $M$ satisfies ACC on annihilators. This proves $\textit{2}$. 

$\textit{2}\Rightarrow\textit{3}:$ Let $N\leq{M}$. Suppose that $N\leq_e{M}$. Since $M$ is semiprime with uniform dimension and satisfies ACC on annihilators, then $M$ is essentially compressible by Proposition  3.13 of \cite{K}. Now, if $f:M\rightarrow{N}$ is a monomorphism then $N\leq_e{M}$ by lemma \ref{132}.

$\textit{3}\Rightarrow\textit{1}:$ It follows from Lemma \ref{107} and Proposition  \ref{125}.
\end{proof}

\newtheorem{20a}[112]{Remark}
\begin{20a}
\normalfont
Notice that Theorem \ref{110} is a generalization of Goldie's Theorem. See \cite{L} Theorem 11.13.

In Proposition  3.13 of \cite{K}, $M$ is a generator of $\sigma[M]$, but by Lemma \ref{148} this hypothesis is not necessary.
\end{20a}

\newtheorem{39}[112]{Corollary}
\begin{39}\label{129}
Let $M$ be projective in $\sigma[M]$ and semiprime. Then, $M$ has finite uniform dimension and enough monoforms if and only if $M$ is a Goldie module.
\end{39}

\begin{proof}
$\Rightarrow:$ Since $M$ is semiprime with finite uniform dimension and enough monoforms, then $M$ is non $M$-singular by Proposition  3.8 of \cite{K}. By Theorem \ref{110}, $M$ is a Goldie module. 

$\Leftarrow:$ If $M$ is a Goldie module, $M$ has finite uniform dimension and by Theorem \ref{110} $M$ is non $M$-singular. Hence the uniform submodules of $M$ are monoform. Since $M$ has finite uniform dimension every submodule of $M$ contains a uniform, hence every submodule contains a monoform.
\end{proof}

For the definition of $M$-Gabriel dimension see \cite{P} section 4.

\newtheorem{45}[112]{Corollary}
\begin{45}\label{135}
Let $M$ be projective in $\sigma[M]$ with finite uniform dimension. If $M$ is a semiprime module and has $M$-Gabriel dimension, then $M$ is a Goldie module. 
\end{45}

\begin{proof}
Let $N\leq{M}$. Since $M$ has $M$-Gabriel dimension, by Lemma 4.2 of \cite{P}, $N$ contains a cocritical submodule $L$. Then $L$ is monoform. By Proposition  \ref{129} $M$ is a Goldie module.
\end{proof}

\newtheorem{40}[112]{Corollary}
\begin{40}\label{130}
Let $M$ be projective in $\sigma[M]$ and semiprime with Krull dimension. Then $M$ is a semiprime Goldie module.
\end{40}

\begin{proof}
Since $M$ has Krull dimension, $M$ has finite uniform dimension and enough monoforms. By Proposition  \ref{129} $M$ is a Goldie module.
\end{proof}

\newtheorem{50}[112]{Proposition}
\begin{50}\label{140}
Suppose that $M$ is progenerator of $\sigma[M]$. Let $N\in\sigma[M]$, then
\[\mathcal{Z}(N)=\sum\{f(M)|f:M\rightarrow{N}\;ker(f)\leq_e{M}\}.\]
\end{50}
\begin{proof}
By definition of $M$-singular module, it is clear that $\sum\{f(M)|f:M\rightarrow{N}\;ker(f)\leq_e{M}\}\leq{\mathcal{Z}(N)}$. Now, let $n\in{\mathcal{Z}(N)}$ and consider $Rn\leq{\mathcal{Z}(N)}$. Since $Rn\in\sigma[M]$ there exists a natural number $t$ and an epimorphism $\rho:M^t\rightarrow{Rn}$. Suppose that $(m_1,..,m_t)$ is such that $\rho(m_1,...,m_t)=n$. If $j_i:M\rightarrow{M^t}$ are the inclusions $(i=1,...,t)$, then by Proposition  \ref{142} $Ker(\rho\circ{j_i})\leq_e{M}$. Thus, $n=\sum_{i=1}^{t}{\rho\circ{j_i}(m_i)}\in\sum\{f(M)|f:M\rightarrow{N}\;ker(f)\leq_e{M}\}$.
\end{proof}

\newtheorem{50k}[112]{Remark}
\begin{50k}\label{140k} 
\normalfont
Let $M\in{R-Mod}$ and consider $\tau_g\in{M-tors}$, where $\tau_g=\xi(\{S\in\sigma[M]|S\;is\;M-singular\})$. If $M\in\mathbb{F}_{\tau_g}$, by \cite{W} Proposition. 10.2, we have that $\chi(M)=\tau_g$. Let $t_{\tau_g}$ be the preradical associated to $\tau_g$. Then
\[t_{\tau_g}(N)=\sum\{S\leq{N}|S\in\mathbb{T}_{\tau_g}\}=\sum\{S\leq{N}|S\;is\;M-singular\}=\mathcal{Z}(N).\]  
\end{50k}

\newtheorem{51}[112]{Proposition}
\begin{51}\label{141}
Suppose $M$ is progenerator of $\sigma[M]$. If $M$ is semiprime Goldie, then 
\[\mathcal{Z}(N)=\sum{f(M)}\]
where the sum is over the $f:M\rightarrow{N}$ such that there exists $\alpha\in{End_R(M)}$ monomorphism with $\alpha(M)\leq_eM$ and $f\alpha=0$.
\end{51}
\begin{proof}
Let $N\in\sigma[M]$. By Proposition  \ref{140}
\[\mathcal{Z}(N)=\sum\{f(M)|f:M\rightarrow{N}\;ker(f)\leq_e{M}\}.\]
If $f:M\rightarrow{N}$ with $Ker(f)\leq_e{M}$, by Theorem \ref{110} there exists a monomorphism $\alpha:M\rightarrow{Ker(f)}$. We have that $f\alpha=0$ and by Lemma \ref{132} $\alpha(M)\leq_e(M)$. 

Let $f:M\rightarrow{N}$ such that there exists $\alpha:M\rightarrow{M}$ $f\alpha=0$ and $\alpha(M)\leq_e(M)$. Then $\alpha(M)\leq{Ker(f)}$. Therefore $Ker(f)\leq_e(M)$.
\end{proof}

\newtheorem{41a}[112]{Remark}
\begin{41a}
\normalfont
Let $R$ be a ring such that $R$-Mod has an infinite set of non-isomorphic simples modules. Consider $M=\bigoplus_I{S_i}$, $I$ an infinite set, such that $S_i$ is a simple module for all $i\in{I}$ and with $Si\ncong{S_j}$ if $i\neq{j}$. This module does not have finite uniform dimension and, in $M$-tors, $\tau_g=\chi$. Then, if $N\in\sigma[M]$
\[t_{\tau_g}(N)=\mathcal{Z}(N)=\sum{f(M)}\]
where the sum is over the $f:M\rightarrow{N}$ such that there exists $\alpha\in{End_R(M)}$ monomorphism with $\alpha(M)\leq_eM$ and $f\alpha=0$.

This example shows that the converse of the last Proposition  is not true in general.
\end{41a}

Following \cite{A}
\newtheorem{56}[112]{Definition}
\begin{56}
A module $M$ is \emph{weakly compressible} if for any nonzero submodule $N$ of $M$, there exists $f:M\rightarrow{N}$ such that $f\circ{f}\neq{0}$.
\end{56}

\newtheorem{57}[112]{Remark}
\begin{57}\label{146}
\normalfont
Notice that if $M$ is weakly compressible then $M$ is a semiprime module. The converse hold if $M$ is projective in $\sigma[M]$
\end{57}

Next definition was taken from \cite{H}
\newtheorem{58}[112]{Definition}
\begin{58}
A module $M$ is a \emph{semiprojective} module if $I=Hom(M,IM)$ for any cyclic right ideal $I$ of $End_R(M)$
\end{58}
For other characterizations see \cite{Foun}.

\newtheorem{19}[112]{Proposition}
\begin{19}\label{109}
Let $M$ be projective in $\sigma[M]$ and retractable. Then, $S:=End_R(M)$ is semiprime if and only if $M$ is semiprime.
\end{19}
\begin{proof}
$\Rightarrow:$ Corollary \ref{108}. 

$\Leftarrow:$ If $M$ is semiprime, since $M$ is projective in $\sigma[M]$ then $M$ is weakly compressible and semiprojective. Then, by [\cite{H}. Theorem 2.6 (b)] $S$ is semiprime.
\end{proof} 

\newtheorem{58a}[112]{Lemma}
\begin{58a}\label{147}
Let $M$ be projective in $\sigma[M]$ and retractable. $M$ is non $M$-singular if and only if $Hom_R(M/N,M)=0$ for all $N\leq_e{M}$.
\end{58a}

\begin{proof}
$\Rightarrow$: If $N\leq_e{M}$ then $M/N$ is $M$-singular, then $Hom_R(M/N,M)=0$. 

$\Leftarrow$: Suppose $\mathcal{Z}(M)\neq{0}$. Since $M$ is retractable there exists $0\neq f:M\to \mathcal{Z}(M)$. By Proposition \ref{142} $Ker(f)\leq_eM$, so there exists a non zero morphism form $M/Ker(f)\to M$.
\end{proof}

For a retractable $R$-module $M$, Theorem 11.6 of \cite{W} gives necessary and sufficient conditions in order to $T:=End_R(\widehat{M})$ being semisimple, left artinian, and being the classical left quotient ring of $S=End_R(M)$. Also, in (\cite{H}, Corollary 2.7) the authors give necessary and sufficient conditions for a semiprojective module $M$ to $S$ being a semiprime right Goldie ring. We give an extension of these results.

\newtheorem{54}[112]{Theorem}
\begin{54}\label{144}
Let $M$ be projective in $\sigma[M]$, $S=End_R(M)$ and $T=End_R(\widehat{M})$. The following conditions are equivalent:
\begin{enumerate}
	\item $M$ is a semiprime Goldie module.
	\item $T$ is semisimple right artinian and is the classical right quotient ring of $S$.
	\item $S$ is a semiprime right Goldie ring.
	\item $M$ is weakly compressible with finite uniform dimension, and for all $N\leq_e{M}$, $Hom_R(M/N,M)=0$ .
\end{enumerate}
\end{54}

\begin{proof}
$\textit{1}\Rightarrow\textit{2}:$  By Proposition  \ref{109}, $S$ is a semiprime ring. Since $M$ is a Goldie module, then $M$ is non $M$-singular with finite  uniform dimension, hence by \cite{W} Proposition 11.6, $T$ is right semisimple and is the classical right quotient ring of $S$.

$\textit{2}\Rightarrow\textit{3}:$ By \cite{L} Theorem 11.13, $S$ is a semiprime right Goldie ring . 

$\textit{3}\Rightarrow\textit{4}:$ By \cite{H} Corollary 2.7.

$\textit{4}\Rightarrow\textit{1}:$ Since $M$ is weakly compressible then $M$ is semiprime. By Lemma \ref{147} $M$ is non $M$-singular. Thus, by Theorem \ref{110} $M$ is a Goldie module.
\end{proof}

\newtheorem{55}[112]{Corollary}
\begin{55}\label{145}
Let $M$ be projective in $\sigma[M]$, $S=End_R(M)$ and $T=End_R(\widehat{M})$. The following conditions are equivalent:
\begin{enumerate}
	\item $M$ is a prime Goldie module.
	\item $T$ is simple right artinian and is the classical right quotient ring of $S$.
	\item $S$ is a prime right Goldie ring.
	\item Given nonzero submodules $N$, $K$ of $M$ there exists a morphism $f:M\rightarrow{N}$ such that $K\nsubseteq{Ker(f)}$. $M$ has finite uniform dimension and for all $N\leq_e{M}$, $Hom(M/N,M)=0$.
\end{enumerate}
\end{55}

\begin{proof}
$\textit{1}\Rightarrow\textit{2}:$ By Proposition \ref{144}, $S$ is a semiprime ring and $T$ is right semisimple and the classical  right quotient ring of $S$. Let $0\neq{I}\leq{T}$ be an ideal. Since $T$ is semisimple, there exits an ideal $J\leq{T}$ such that $T=I\oplus{J}$. Put $M_1=I\widehat{M}$ and $M_2=J\widehat{M}$. Then $M_1$ and $M_2$ are fully invariant submodules of $\widehat{M}$ and $M_1\cap{M_2}=0$ because $I\cap{J}=0$. Consider $M_1\cap{M}$ and $M_2\cap{M}$. If $f\in{S}$, then there exists $\hat{f}\in{T}$ such that $f=\hat{f}|_M$. Let $x\in{M_1\cap{M}}$. Then $f(x)=\hat{f}(x)\in{M_1\cap{M}}$ since $M_1$ is a fully invariant submodule of $\widehat{M}$. Thus $M_1\cap{M}$ is a fully invariant submodule of $M$. In the same way, $M_2\cap{M}$ is fully invariant in $M$. Since $(M_1\cap{M})\cap{(M_2\cap{M})}=0$, then $(M_1\cap{M})_M{(M_2\cap{M})}=0$. Hence $M_1\cap{M}=0$ or $M_2\cap{M}=0$ because $M$ is prime. On the other hand, $M\leq_e\widehat{M}$ and so $M_1=0$ or $M_2=0$. Since $0\neq{I}$, then $M_2=0$. Thus $J=0$, and it follows that $T$ is a simple ring. 

$\textit{2}\Rightarrow\textit{3}:$ By \cite{L} Corollary 11.16, $S$ is a prime right Goldie ring. 

$\textit{3}\Rightarrow\textit{4}:$ Let $N$, $K$ be nonzero submodules of $M$, if $K\subseteq{Ker(f)}$ for all $f:M\rightarrow{N}$ then $0=Hom_R(M,N)Hom(M,K)\leq{S}$. Then $Hom_R(M,N)=0$ or $Hom_R(M,K)=0$. By retractability, $N=0$ or $K=0$, a contradiction.

$\textit{4}\Rightarrow\textit{1}$ It is clear.
\end{proof}

\newtheorem{49a}[112]{Remark}
\begin{49a}\label{sing}
\normalfont
Suppose that $M$ and $N$ are $R$-modules such that $\sigma[N]\subseteq\sigma[M]$. If $N$ is non $M$-singular, then $N$ is non $N$-singular. This is because if there exists an exact sequence $0\rightarrow{L}\rightarrow{K}\rightarrow{N}\rightarrow{0}$ in $\sigma[N]$ such that $L\leq_e{N}$, then this sequence is in $\sigma[M]$ which implies that $N$ is $M$-singular, a contradiction.
\end{49a}

\newtheorem{49}[112]{Proposition}
\begin{49}\label{139}
Let $M$ be projective in $\sigma[M]$ and semiprime with finitely many minimal prime submodules $P_1,...,P_t$. Suppose every quotient $M/P_{i}$ $(1\leq{i}\leq{t})$ has finite uniform dimension. Then $M$ is a Goldie module if and only if each $M/P_i$ is a Goldie module.
\end{49}

\begin{proof}
$\Rightarrow:$ Suppose $M$ is a Goldie module and $P_i$ is a minimal prime submodule of $M$. By hypothesis, each $M/P_i$ has finite uniform dimension. Notice that by proposition \ref{116a}
\[P_i\subseteq{Ann_M(P_1\cap...\cap{P_{i-1}}\cap{P_{i+1}}\cap...\cap{P_n})}\]
Since $M$ has finite uniform dimension there exist a uniform submodule $U_i$ of $P_1\cap...\cap{P_{i-1}}\cap{P_{i+1}}\cap...\cap{P_n}$. So $P_i\subseteq{Ann_M(U_i)}$. By Propositions \ref{115} and \ref{115b}, $P_i=Ann_M(U_i)$. Then, there exists a monomorphism $M/P_i\rightarrow{{U_i}^X}$ and since $U_i$ is non $M$-singular, then $M/P_i$ is non $M$-singular. Thus $M/P_i$ is non $({M}/{P_i})$-singular by Remark \ref{sing}. Since $M/P_i$ is a prime module, by Theorem \ref{110} $M/P_i$ is a Goldie module. 

$\Leftarrow:$ By Corollary \ref{116a} there exists a monomorphism $M\rightarrow\bigoplus_{i=1}^{t}{M/P_i}$. Since each $M/P_i$ has finite uniform dimension then $M$ has finite uniform dimension. 

Let $0\neq{N}$ be a submodule of $M$. Since there exists a monomorphism $M\rightarrow\bigoplus{M/P_i}$ then there exists $1\leq{i}\leq{t}$ and submodules $0\neq{K}\leq{M/P_i}$ and $0\neq{N'}\leq{N}$ such that $K\cong{N'}$. We have that $M/P_i$ is a Goldie module, thus it has enough monoforms. Hence $N'$ has a monoform submodule, that is $M$ has enough monoforms, and so by Corollary \ref{129} $M$ is Goldie module.
\end{proof}

\newtheorem{42a}[112]{Remark}
\begin{42a}\label{132a}
\normalfont
Notice that if $M$ is a semiprime Goldie module then $M$ has finitely many minimal prime submodules by Proposition  \ref{117}. So in the proof $\Rightarrow:$ of Proposition  \ref{139} this hypothesis is not used.
\end{42a}

\newtheorem{41}[112]{Definition}
\begin{41}\label{131}
\normalfont
Let $M\in{R-Mod}$ and $N\leq{M}$. We say $N$ is a \emph{regular submodule}  if there exists a monomorphism $M\rightarrow{N}$. Denote 
\[\textit{Reg}(M):=\{N\leq{M}|N\;regular\;submodule\}\]
\end{41}

\newtheorem{b1}[112]{Remark} 
\begin{b1}\label{bb}
\normalfont
There exists modules with regular submodules which are nonessential. For example, a pure infinite module, see \cite{MM}.
\end{b1}

\newtheorem{43}[112]{Proposition}
\begin{43}\label{133}
Let $M$ be projective in $\sigma[M]$ and a semiprime Goldie module. Then, $N$ is a regular submodule of $M$ if and only if $N$ is essential in $M$.
\end{43}

\begin{proof}
Since $M$ is Goldie, every regular submodule is essential by Lemma \ref{132}. Now, let $N\leq_eM$. By Theorem \ref{110}, $N$ is a regular submodule.
\end{proof} 

If $K\in\sigma[M]$, we say that $K$ is $\textit{Reg}(M)$-injective if any morphism $f:N\to K$ with $N\in\textit{Reg}(M)$ can be extended to a endomorphism of $M$.

\newtheorem{44}[112]{Corollary}
\begin{44}\label{134}
Let $M$ be projective in $\sigma[M]$ and a semiprime Goldie module. Let $K\in\sigma[M]$. If $K$ is $\textit{Reg}(M)$-injective, then $K$ is $M$-injective.
\end{44}

\section{Duo Modules}

Following \cite{DM}
\newtheorem{114}{Definition}[section]
\begin{114}\label{118}
\normalfont
Let $M\in{R-Mod}$. $M$ is a \emph{duo module}  if every submodule of $M$ is fully invariant in $M$.
\end{114}

Examples:
\begin{enumerate}
	\item If $_RS$ is a simple module then, $S$ is a duo module.
	\item If $_RM=\bigoplus_I{S_i}$ with $S_i$ simple and $S_i$ not isomorphic to $S_j$ $i\neq{j}$ then $M$ is a duo module.
	\item An $R$-module $M$ is called a multiplication module if every $N\leq{M}$ is of the form $IM=N$ for some ideal $I$ of $R$. These modules are examples of duo modules. See \cite{T}
	\item Consider the example in Remark \ref{aa} that was taken from \cite{P}. In that paper it is proved that $M/K\cong{S}\cong{M/L}\cong{M/N}$, hence $L$, $K$ and $N$ are maximal submodules of $M$. It follows that $K/S$, $L/S$ and $N/S$ are maximal submodules of $M/S$. Moreover, since $K\cap{L}=S=K\cap{N}=N\cap{L}$, then $M/S=K/S\oplus{L/S}$. Thus
	\[K/S\cong\frac{M/S}{L/S}\cong{M/L}\cong{S}\]
	This implies that $K/S$ is simple, and analogously $L/S$ and $N/S$ are simple. 
	
	Let $0\neq{T}<{M}$. Since $S\leq_e{M}$, then $S\leq{T}$. If $T=S$, then $T$ is fully invariant. Suppose that $T\neq{S}$ and $T\notin\{K,L,N\}$. We have that $S\leq{T\cap{K}}\leq{K}$. Moreover, since $K/S$ is simple, then $T\cap{K}=S$ or $T\cap{K}=K$. If $T\cap{K}=K$ then $K\leq{T}<M$; but $K$ is maximal, then $K=T$, a contradiction. Thus, $T\cap{K}=S$. Analogously $T\cap{L}=S=T\cap{N}$. 
	
	Let $0\neq{x}\in{M}$. If $ann_R(x)=0$, there exists a monomorphism $R\rightarrow{M}$ and thus $E(R)=M$, a contradiction, because $E(R)\cong{M}\oplus{M}$ (see \cite{P}, Example 1.12) and $M$ is a indecomposable injective module. Thus, $ann_R(x)\neq{0}$ for all $0\neq{x}\in{M}$. 
	
	Let $0\neq{x}\in{T}$. Since $ann_R(x)\neq{0}$, then $ann_R(x)\in\{I,J_1,J_2,J_3\}$. By Theorem 2.13 of \cite{BP} we have that:
\begin{itemize}
\item	If $ann_R(x)=I$ then $x\in{S}$
	
\item	If $ann_R(x)=J_1$ then $x\in{K\cap{T}}=S$ 
	
\item	If $ann_R(x)=J_2$ then $x\in{L\cap{T}}=S$ 
	
\item	If $ann_R(x)=J_3$ then $x\in{N\cap{T}}=S$ 
\end{itemize}

	Therefore $T\leq{S}$, a contradiction. Thus, all submodules of $M$ are fully invariant. 
\end{enumerate}

\newtheorem{46a}[114]{Remark}
\begin{46a}
\normalfont
	In \cite{DM} the authors state that they did not know an example of a duo module $M$ and a submodule $N$ such that $M/N$ is not a duo module. In this example, $M$ is a duo module, but $M/S\cong{S}\oplus{S}$ is not a duo module.  
\end{46a}

\newtheorem{46}[114]{Proposition}
\begin{46}\label{136}
$M$ is a duo module as $R$-module and it generates all its submodules if and only if $M$ is a multiplication module as $End_R(M)$-module.
\end{46}

\begin{proof}
$\Rightarrow:$ Let $S=End_R(M)$ and let $N$ be a submodule of $M$. Since $M$ is a duo module, $N$ is fully invariant, thus $Hom_R(M,N)$ is an ideal of $S$. Since $M$ generates all its submodules, then $N=tr^M(N)=Hom_R(M,N)M$. Thus, $M$ is a multiplication module as $End_R(M)$-module. 

$\Leftarrow:$ It is clear.
\end{proof}

\newtheorem{29}[114]{Proposition}
\begin{29}\label{119}
Let $M$ be projective in $\sigma[M]$. Suppose that $M$ is a semiprime and non $M$-singular duo module. Then, for every subset $X\subseteq{End_R(M)}$ we have that:
\[Ann_M(Ann_M(\bigcap_X{Ker(f)}))=\bigcap_X{Ker(f)}\]
\end{29}

\begin{proof}
Since $M$ is a duo module, by Proposition  \ref{114}, $Ann_M(\bigcap_X{Ker(f)})$ is the unique pseudocomplement of $\bigcap_X{Ker(f)}$. Then
\[\bigcap_X{Ker(f)}\leq_eAnn_M(Ann_M(\bigcap_X{Ker(f)})).\]
Since $M$ is non $M$-singular, $\bigcap_X{Ker(f)}$ has no essential extensions in $M$ by Lemma 3.5 of \cite{K}. Thus, we have the equality.
\end{proof}

\newtheorem{30}[114]{Proposition}
\begin{30}\label{120}
Let $M$ projective in $\sigma[M]$. Suppose that $M$ is a semiprime and non $M$-singular duo module. The following conditions are equivalent:
\begin{enumerate}
	\item $M$ has finite uniform dimension.
	\item $M$ has a finite number of minimal prime submodules.
	\item The number of annihilators in $M$ is finite.
	\item $M$ satisfies the ACC on annihilators.
	\item $M$ satisfies the ACC on pseudocomplements.  
\end{enumerate}
\end{30}

\begin{proof}
$\textit{1}\Rightarrow\textit{2}\Rightarrow\textit{3}:$ Are true by Proposition  \ref{117}. 

$\textit{3}\Rightarrow\textit{4}:$ By Proposition  \ref{119}. 

$\textit{4}\Rightarrow\textit{5}:$ By Proposition  \ref{114}.

$\textit{5}\Rightarrow\textit{1}:$ By \cite{L} Proposition 6.30.
\end{proof}

\newtheorem{47}[114]{Proposition}
\begin{47}\label{137}
Let $M$ be projective in $\sigma[M]$. Suppose $M$ is a prime duo module with finite uniform dimension. Then, $Udim(M)=1$
\end{47}

\begin{proof}
Since $M$ is prime, $0$ is the unique minimal prime submodule of $M$. By Proposition  \ref{117}, there exists a uniform submodule $U$ of $M$ such that $0=Ann_M(U)$. By Proposition  \ref{114}, $U\leq_eM$. Thus, $Udim(M)=1$.
\end{proof}

\newtheorem{48}[114]{Theorem}
\begin{48}\label{138}
Let $M$ be projective in $\sigma[M]$. If $M$ is a semiprime duo module, then the following conditions are equivalent: 
\begin{enumerate}
	\item $M$ is a prime Goldie module.
	\item $\widehat{M}$ is indecomposable and $M$ is non $M$-singular.
	\item $M$ is uniform and non $M$-singular.
\end{enumerate}
\end{48}

\begin{proof}
$\textit{1}\Rightarrow\textit{2}:$ Since $M$ is a prime module, by Proposition  \ref{137}, $Udim(M)=1$ and then $\widehat{M}$ is indecomposable. Since $M$ is a Goldie module, by Theorem \ref{110} $M$ is non $M$-singular. 

$\textit{2}\Rightarrow\textit{3}:$ Let $0\neq{K}\leq{M}$. Then, there exists $L\leq{M}$ such that $K\oplus{L}\leq_eM$. Hence, $\widehat{K}\oplus\widehat{L}=\hat{M}$, but since $\widehat{M}$ is indecomposable, then $L=0$. Thus $K\leq_e{M}$. 

$\textit{3}\Rightarrow\textit{1}:$ Let $K$ and $0\neq{L}$ be submodules of $M$ such that $K_ML=0$. Then, $K\leq{Ann_M(L)}$, and thus $K\cap{L}=0$ by Proposition  \ref{114}. Since $M$ is uniform, $K=0$. Thus, $M$ is prime and by Theorem \ref{110} $M$ is Goldie.
\end{proof}

\end{document}